\documentclass[10pt]{amsart}

\usepackage{amsthm,amsmath,amssymb,amscd,amsfonts,latexsym}

\theoremstyle{plain}
\newtheorem{theorem}{Theorem}[section]

\newtheorem{maintheorem}[theorem]{Main Theorem}
\newtheorem*{maintheorem12}{Main Theorem 1.2}
\newtheorem*{theorem13}{Theorem 1.3}

\newtheorem{proposition}[theorem]{Proposition}
\newtheorem{criterion}[theorem]{Criterion}
\newtheorem{corollary}[theorem]{Corollary}
\newtheorem{def-thm}[theorem]{Definition-Theorem}
\newtheorem{lemma}[theorem]{Lemma}

\theoremstyle{definition}
\newtheorem{definition}[theorem]{Definition}

\newtheorem{remark}[theorem]{Remark}

\newtheorem{facts}[theorem]{Facts}

\newtheorem*{acknowledgement}{Acknowledgement}

\newcommand{\PP}{\mathbb{P}}

\newcommand{\RR}{\mathbb{R}}

\newcommand{\NN}{\mathbb{N}}
\newcommand{\ZZ}{\mathbb{Z}}
\newcommand{\CC}{\mathbb{C}}

\newcommand{\OO}{{\mathcal O}}

\newcommand{\II}{{\mathcal I}}

\newcommand{\ke}{{K\"ahler-Einstein\ }}

\DeclareMathOperator{\mult}{mult}

\DeclareMathOperator{\Pic}{Pic}

\DeclareMathOperator{\Aut}{Aut}
\DeclareMathOperator{\Cr}{Cr}
\DeclareMathOperator{\id}{id}
\DeclareMathOperator{\Ricci}{Ricci}
\DeclareMathOperator{\rank}{rank}
\DeclareMathOperator{\PGL}{PGL}
\DeclareMathOperator{\GL}{GL}

\begin{document}

\title[Existence of K\"ahler-Einstein metrics and multiplier ideal sheaves]
{Existence of K\"ahler-Einstein metrics and multiplier ideal sheaves on del Pezzo surfaces}
\begin{abstract} 
We apply Nadel's method of multiplier ideal sheaves to show that every complex del Pezzo surface of degree at most six whose automorphism group acts without fixed points has a K\"ahler-Einstein metric. In particular, all del Pezzo surfaces of degree $4,5,$ or $6$ and certain special del Pezzo surfaces of lower degree are shown to have a K\"ahler-Einstein metric. This result is not new, but the proofs given in the present paper are less involved than earlier ones by Siu, Tian and Tian-Yau.
\end{abstract}

\author[G. Heier]{Gordon Heier}
\address{Department of Mathematics\\ Christmas-Saucon Hall\\ Lehigh University\\ \ \ 14 E.~Packer Avenue \\
Bethlehem, PA 18015\\ USA}

\email{heier@lehigh.edu}


\maketitle

\section{Introduction} 
Let $X$ be a compact complex manifold. It is said to be a {\it \ke manifold} if there exists a K\"ahler form $\omega=\frac{i}{2}\sum_{i,j} g_{i\bar{j}} dz_i\wedge d\bar{z}_{j}$ such that $\omega$ and its Ricci curvature form are proportional, ie
\begin{equation*}
\Ricci(\omega)=-i\partial\bar\partial \log\det(g_{i\bar{j}})=\lambda \omega
\end{equation*}
for some $\lambda \in \RR$. The hermitian metric $\sum_{i,j} g_{i\bar{j}} dz_i\otimes d\bar{z}_{j}$ corresponding to such an $\omega$ is called a {\it \ke metric}. Since the anticanonical class $c_1(-K_X)$ (which is also referred to as the first Chern class $c_1(X)$) is known to always contain the form $\frac{1}{2\pi}\Ricci(\omega)$ by a result of Chern, \ke metrics can only exist on manifolds with definite or vanishing anticanonical class.\par
In his ground-breaking work \cite{Yau_PNAS,Yau_Comm_PAM}, Yau proved Calabi's conjecture, which says that on a compact K\"ahler manifold, for any given K\"ahler class, every form representing $2\pi c_1(X)$ can be obtained as the Ricci curvature form of a unique K\"ahler form in the given K\"ahler class. As a corollary, on a K\"ahler manifold with vanishing first Chern class, to every K\"ahler class there corresponds a unique Ricci flat \ke metric.\par
On compact complex manifolds with negative anticanonical line bundle, the existence of K\"ahler-Einstein metrics has been established independently by Aubin \cite{Aubin_CR,Aubin_Bull} and Yau (op.~cit.). In the case of positive anticanonical line bundle, a K\"ahler-Einstein metric may or may not exist. Different obstructions have been formulated by Matsushima \cite{Matsushima}, Futaki \cite{Futaki} and others (see \cite{Bourguignon} for a nice survey).\par
In the mid 1980's, Siu \cite{Siu_ke_symm}, Tian \cite{Tian_1987} and Tian-Yau \cite{Tian_Yau} provided sufficient conditions for the existence of K\"ahler-Einstein metrics in the case of positive anticanonical bundles, settling certain cases (see below). However, these results were not strong enough to completely clarify the question of existence of K\"ahler-Einstein metrics on del Pezzo surfaces. Recall that del Pezzo surfaces are two-dimensional compact complex manifolds with positive anticanonical line bundle. The classification of del Pezzo surfaces tells us that they are isomorphic either to $\PP^1\times\PP^1$,  $\PP^2$, or $\PP^2$ blown up in $r\in\{1,\ldots,8\}$ general points (see Section \ref{delpezzoclass}).\par
For $\PP^1\times\PP^1$ and  $\PP^2$, the existence of a \ke metric follows directly from the fact that the Fubini-Study metric on $\PP^1$ and $\PP^2$ is K\"ahler-Einstein. For $\PP^2$ blown up at one or two points, one can show that the so-called Futaki-invariant does not vanish and therefore obstructs the existence of a \ke metric in these cases (see \cite[Examples 3.10, 3.11]{Tian_book} for a nice discussion). Alternatively, and more in the spirit of this paper, it suffices to observe that the automorphism groups are non-reductive Lie groups, which precludes the existence of \ke metrics according to Matsushima. For the unique del Pezzo surface coming from the blow up of three points, the existence of a \ke metric was shown by \cite{Siu_ke_symm} and, independently, by the combined efforts of \cite{Tian_1987} and \cite{Tian_Yau}. For the unique del Pezzo surface coming from the blow up of four points, the existence of a \ke metric was shown by \cite{Tian_1987} and \cite{Tian_Yau}. For the general case of del Pezzo surfaces coming from the blow up of  $r\in\{5,\ldots,8\}$ points, Tian's paper \cite{Tian_Inv_Math} finally gave a proof of the fact that they do carry \ke metrics. However, the proof certainly is not short in length.
\begin{remark}
In the paper \cite{Tian_Yau}, the authors apparently overlooked the fact that {\it every} del Pezzo surface coming from the blow up of five points in $\PP^2$ can be obtained as the intersection of quadrics in $\PP^4$ of the form
\begin{equation*}
\sum_{i=0}^4X_i^2=\sum_{i=0}^4 a_iX_i^2=0,
\end{equation*}
where $a_i\not=a_j$ for $i\not= j$. This fact follows eg from \cite[p. 551]{GH} together with the Normal Form Lemma for Simple Pencils of Quadrics \cite[Lemma 22.42]{Harris_First_Course}. Therefore, the arguments in \cite[pp. 188-189]{Tian_Yau} actually settle the case $r=5$ completely, and the below-mentioned dichotomy is already present in the work of Tian and Yau.
\end{remark}
A different approach to the problem of existence of \ke metrics was taken by Nadel in his seminal paper \cite{Nadel_Annals}. In that paper, Nadel introduces the notion of a multiplier ideal subsheaf of the sheaf of holomorphic functions, partly motivated by the successful use of a similar concept in Kohn's paper \cite{Kohn_Acta} on boundary regularity for the complex Neumann problem on weakly pseudoconvex domains of finite type. Nadel proves that the non-existence of certain multiplier ideal sheaves is a sufficient condition for the existence of \ke metrics on a given Fano manifold.\par
The purpose of this paper is to answer the following open question: To what extent can the technique of multiplier ideal sheaves be used to prove the existence of K\"ahler-Einstein metrics on del Pezzo surfaces? It turns out that there is not a uniform answer to this question due to a dichotomy between the case of at most five points blown up and the case of at least six points blown up. The precise main theorem we (re-)prove is the following.
\begin{maintheorem} [\cite{Siu_ke_symm}, \cite{Tian_1987}, \cite{Tian_Yau}, \cite{Tian_Inv_Math}]
Let $X$ be a del Pezzo surface obtained by blowing up $\PP^2$ in $3,4,$ or $5$ points. Then $X$ carries a \ke metric.
\end{maintheorem}
Naturally, the reader may wonder what happens for six or more points blown up. Due to the nature of our proof based on automorphism groups, we find it impossible to prove the analogous statement for $6,7,$ or $8$ points using our method. However, it turns out that there is a previous result of Nadel that addresses precisely the case of $6,7,$ or $8$ points blown up:
\begin{theorem}[\cite{Nadel_Wuppertal}] \label{Nadel_low_degree}
Let $X$ be a del Pezzo surface obtained by blowing up $\PP^2$ in $6,7,$ or $8$ points. Let the automorphism group of $X$ act without fixed points on $X$. Then $X$ carries a \ke metric.
\end{theorem}
Note that Nadel's result is complementary to ours in the sense that it requires six or more points in an essential way to work. For the sake of completeness, we briefly reproduce Nadel's proof of Theorem \ref{Nadel_low_degree} in Section \ref{678}, filling in some minor missing details in the proof of Proposition \ref{no_one_dim_mis}. The question regarding the applicability of Theorem \ref{Nadel_low_degree} (eg to the cubic Fermat hypersurface in $\PP^3$) is discussed at the very end of this paper. 

\begin{acknowledgement} The author would like to thank I. Dolgachev for helpful comments concerning the automorphism groups of del Pezzo surfaces of degree at most three.
\end{acknowledgement}

\section{Classification and basic properties of del Pezzo surfaces}\label{delpezzoclass}
\begin{definition}
A {\it del Pezzo surface} is a two-dimensional compact complex manifold $X$ whose anti-canonical line bundle $-K_X$ is ample. We call the self-intersection number $(-K_X)^2=K_X^2$ the {\it degree} of $X$. We will denote the degree also by $d_X$.
\end{definition}
We now gather some important facts about del Pezzo surfaces, resulting in the standard classification (see \cite{Manin, Demazure, Hartshorne}).

\begin{facts}
For every del Pezzo surface $X$, the Picard group $\Pic X$ satisfies
\begin{equation*}
\rank \Pic X + d_X=10.
\end{equation*}
In particular, $d_X\leq 9$.\par
If $d_X=9$, then $X$ is isomorphic to $\PP^2$.\par
If $d_X=8$, then $X$ is isomorphic either to $\PP^1\times\PP^1$ or to $\tilde \PP^2$, ie $\PP^2$ blown up at one point.\par
If $7\geq d_X\geq 1$, then $X$ is isomorphic to $\PP^2$ blown up at $r=9-d_X$ points which have the following properties:
\begin{enumerate}
\item no three points lie on a line,
\item no six points lie on a conic,
\item no seven points lie on a cubic such that the eighth is a double point of the cubic.
\end{enumerate}
Any set of $r=9-d_X$ points satisfying the above three properties will be said to be in {\it general} position, and, conversely, the result of blowing up $1\leq r \leq 8$ points in general position in $\PP^2$ is a del Pezzo surface. For $1\leq r\leq 4$ general points blown up, there is in each case a unique resulting del Pezzo surface. The reason is that for any two sets of points $P_1,\ldots,P_r$ and $Q_1,\ldots,Q_r$ ($r\leq 4$), with each set in general position, there is an element $A\in \Aut(\PP^2)=\PGL(3,\CC)$ with $A(P_i)=Q_i\ (1\leq i\leq r)$.
\end{facts}
For our understanding of del Pezzo surfaces, the following facts about the anti-canonical line bundle are also very important.
\begin{facts}
Let $1\leq r\leq 8$. Let $X$ be obtained by blowing up $r$ general points. Let $E_i$ denote the exceptional $(-1)$-curves that are the pre-images of the $r$ points. Let $\pi:X\to\PP^2$ denote the blow-up map. Then
\begin{equation*}
K_X=\pi^*K_{\PP^2}+\sum_{i=1}^{r}E_i.
\end{equation*}
This yields
\begin{equation*}
\dim H^0(X,-K_X)=10-r.
\end{equation*}\par
For $1\leq r\leq 6$, the complete linear system $|-K_X|$  gives an embedding into $\PP^{9-r}=\PP^{d_X}$. For $r=7$, it gives a double cover of $\PP^2$. The complete linear system $|-2K_X|$ gives an embedding into $\PP^6$. For $r=8$, $|-K_X|$ has a unique base point, $|-2K_X|$ gives a double cover of a singular quadric surface in $\PP^3$, and $|-3K_X|$ gives an embedding into $\PP^6$.\par

Finally, it turns out that, on every del Pezzo surface of degree at most $7$, the number of $(-1)$-curves exceeds $r$. The reason is that, when blowing up two points in $\PP^2$, the proper transform of the unique line through the two points becomes a $(-1)$-curve as well. When blowing up five points, the unique conic through the five points also becomes a $(-1)$-curve. It is easy to count these $(-1)$-curves: for $r=1,\ldots,8$, their numbers are $1,3,6,10,16,27,56,240$, respectively. Interestingly, for $r=1,\ldots,6$, under the map given by $|-K_X|$, all $(-1)$-curves become lines in projective space. Therefore, they are often referred to as {\it lines} on $X$. 
\end{facts}

\begin{remark}
It seems appropriate to remark that most facts about the automorphism groups of del Pezzo surfaces used below were already known in the 19th century (see \cite{Wiman}). The references to more modern treatments given throughout the text are meant for the reader's convenience, not to apportion credit.\par
\end{remark}

\section{Nadel's method of multiplier ideal sheaves}\label{Nadel_method}
The following is the standard definition of the multiplier ideal sheaf pertaining to a plurisubharmonic function on a complex manifold.
\begin{def-thm}[\cite{Nadel_Annals}]
Let $\varphi$ be a plurisubharmonic function on the complex manifold $X$. Then the {\it multiplier ideal sheaf} $\II(\varphi)$ is the subsheaf of $\OO_X$ defined by
\begin{equation*}
\II(\varphi)(U) =\{f\in \OO_X(U):|f|^2e^{-\varphi}\in L^1_{\text{loc}}(U)\}
\end{equation*}
for every open set $U\subseteq X$. It is a coherent subsheaf.
\end{def-thm}
Multiplier ideal sheaves have turned out to be very useful in algebraic geometry, mainly because of the following vanishing theorem. They are usually defined using the notion of a singular hermitian metric on a line bundle, which in general is a metric $h$ that is given on a small open set $U$ by $h=e^{-\varphi}$, where $\varphi$ is $L^1(U)$. If $\varphi$ is plurisubharmonic for every $U$, the multiplier ideal sheaf $\II(h)$ attached to $h$ is defined by $\II(h)(U)=\II(\varphi)(U)$ if $h=e^{-\varphi}$ on $U$.
\begin{theorem}[Nadel's vanishing theorem]
Let $X$ be a compact complex K\"ahler manifold. Let $L$ be a line bundle on $X$ equipped with a singular hermitian metric such that the curvature current $-\frac{i}{2\pi}\partial\bar\partial \log h$ is positive definite in the sense of currents, ie there is a smooth positive definite $(1,1)$-form $\omega$ and $\varepsilon >0$ such that $-\frac{i}{2\pi}\partial\bar\partial \log h \geq \varepsilon\omega$. Then
\begin{equation*}
H^q(X,(K_X+L)\otimes \II(h))=0 \quad \text{for all } q\geq 1.
\end{equation*}
\end{theorem}
It is well-known that the existence of a \ke metric is equivalent to the solvability of a certain Monge-Amp\`ere equation. If no \ke metric exists, then the continuity method to solve the Monge-Amp\`ere equation must fail. The key idea of \cite{Nadel_Annals} is to capture the failure in the form of an invariant non-trivial multiplier ideal sheaf, as expressed in the following theorem (see also \cite{Demailly_Kollar} for a concise presentation of the details).

\begin{theorem}[Nadel's existence criterion for \ke metrics]\label{Nadel_crit}
Let $X$ be a Fano manifold, ie let $-K_X$ be positive. Assume that $X$ does not have  a \ke metric. Then, for all compact $G\subseteq \Aut(X)$, the line bundle $-K_X$ possesses a $G$-invariant singular hermitian metric $h=h_0e^{-\varphi}$, with $h_0$ a smooth $G$-invariant metric of smooth positive definite curvature $\omega_0$ and $\varphi\in L^1_{\text{loc}}(X)$ $G$-invariant, such that
\begin{enumerate}
\item the curvature current $\Theta_{h}$ of $h$ satisfies
\begin{equation*}
\Theta_{h}=-\frac{i}{2\pi}\partial\bar\partial \log h=\omega_{0}+\frac{i}{2\pi}\partial\bar\partial \varphi \geq 0,
\end{equation*}
\item $\forall \gamma\in ]\frac n {n+1},1[: 0\not =\II(\gamma \varphi)\not = \OO_X$.
\end{enumerate}
The multiplier ideal sheaf $\II(\gamma \varphi)$ is also $G$-invariant. In particular, every element of $G$ maps the zero-set $V(\II(\gamma \varphi))$ to itself.
\end{theorem}
Using the Nadel vanishing theorem, Nadel's criterion yields the following corollaries.
\begin{corollary}
Let $X, G, h_0,\omega_{0},\varphi,h,\Theta_{h},\gamma$ be as in Theorem \ref{Nadel_crit} (Nadel's criterion). Then
\begin{equation*}
H^q(X,\II(\gamma\varphi))=0 \quad \forall q \geq 1.
\end{equation*}
\end{corollary}
\begin{proof}
The proof consists of applying Nadel's vanishing theorem with $L=-K_X$, $h_\gamma=h_0e^{-\gamma\varphi}$. In order to do that, we have to verify that the curvature satisfies the assumption in Nadel's vanishing theorem:
\begin{eqnarray*}
-\frac{i}{2\pi}\partial\bar\partial \log h_\gamma&=& -\frac{i}{2\pi}\partial\bar\partial\log h_0e^{-\gamma\varphi}\\
&=& -(1-\gamma)\frac{i}{2\pi}\partial\bar\partial \log h_0 -\gamma \frac{i}{2\pi}\partial\bar\partial \log (h_0e^{\varphi})\\
&\geq& (1-\gamma)\omega_{0} + \gamma  \Theta_{h}\\
&\geq& (1-\gamma)\omega_0.
\end{eqnarray*}
Now that we have established the required positivity, it can be concluded that 
\begin{equation*}
H^q (X,(K_X-K_X)\otimes \II(\gamma\varphi))=H^q (X,\II(\gamma\varphi))=0 \quad \forall q\geq 1.
\end{equation*}
\end{proof} 
\begin{corollary}\label{van_with_-K}
Let $X, G, h_0,\omega_{0},\varphi,h,\Theta_{h},\gamma$ be as in Theorem \ref{Nadel_crit} (Nadel's criterion). Let ${V_\gamma}=V(\II(\gamma\varphi))\not = \emptyset$. Then 
\begin{equation*}
H^q({V_\gamma},\OO_{V_\gamma})=\begin{cases} 
\CC& \text{if } q=0\\ 
0& \text{if } q \geq 1. 
\end{cases} 
\end{equation*} 
\end{corollary} 
\begin{proof}
Since $-K_X$ is ample, it follows from Kodaira's Vanishing Theorem that
\begin{equation*}
H^q(X,K_X-K_X)=H^q(X,\OO_X)=0\quad \forall q\geq 1.
\end{equation*}
Applying this together with Nadel's Vanishing Theorem to the long exact sequence of 
\begin{equation*}
0\to \II(\gamma\varphi)\to\OO_X\to\OO_{V_\gamma}\to 0
\end{equation*}
yields the result.
\end{proof}
Corollary \ref{van_with_-K} tells us that if $\dim {V_\gamma}=0$, then $V_\gamma$ is an isolated point that is a fixed point of $G$, ie the point constitutes its full orbit under the action of $G$. Moreover, if $\dim {V_\gamma}=1$, it follows that ${V_\gamma}$ is a tree of smooth rational curves (see \cite[Theorem 4.4]{Nadel_Annals}).\par

Now, here is how we will show the existence of a \ke metrics on the del Pezzo surfaces under discussion.\par

First, we consider the case $r\leq 5$.\par

The key observation in this case is that Theorem \ref{Nadel_crit} can be reduced to Criterion \ref{simple_crit} below, based on the following Theorem \cite[Theorem 4.5]{Nadel_Annals}. It is proven by an induction argument.
\begin{theorem}\label{invar_comp}
Let $G\subset \Aut(X)$ be a finite subgroup acting without fixed points on $V_\gamma$. Then some irreducible component of $V_\gamma$ is invariant under $G$. This component is a smooth rational curve.
\end{theorem}
Theorem \ref{invar_comp} allows us to formulate the following simple criterion for the existence of \ke metrics on del Pezzo surfaces.
\begin{criterion}\label{simple_crit}
Let $X$ be a del Pezzo surface. Let $G$ be a finite subgroup of $\Aut(X)$, acting on $X$ without fixed points. Let $G$ act effectively on all $G$-invariant curves on $X$. Let $G$ be not isomorphic to any one of the following finite groups: 
\begin{enumerate}
\item the cyclic group: $\ZZ_n (n\in\NN)$,
\item the dihedral group: $D_{2n}=\ZZ_n\rtimes \ZZ_2 (n\in\NN)$,
\item the alternating group on four letters: $A_4$,
\item the alternating group on five letters: $A_5$,
\item the symmetric group on four letters: $S_4$.
\end{enumerate}
Then $X$ has a \ke metric.
\end{criterion}
\begin{proof}
Let us assume that there exists a one-dimensional $G$-invariant zero-set of a multiplier ideal sheaf. Since $G$ acts effectively on all $G$-invariant curves on $X$, in particular it will act effectively on the $G$-invariant smooth rational curve whose existence is established in Theorem \ref{invar_comp}. However, the finite subgroups of $\Aut(\PP^1)$ are precisely the ones given in the above list (see \cite[Chapter 2]{Grove_Benson}), which yields a contradiction. \par
Since $G$ acts fixed point free, there also are no zero-dimensional zero-sets of multiplier ideal sheaves. Now Theorem \ref{Nadel_crit} (Nadel's criterion) yields the existence of a \ke metric.
\end{proof}

Secondly, we consider the case of $6$ or more points.\par

We need the following relative version of Corollary \ref{van_with_-K}, whose proof commonly is based on the Leray spectral sequence  (see \cite{Nadel_Wuppertal}, \cite{PAGII}). 
\begin{corollary}
In the situation of Theorem \ref{Nadel_crit} (Nadel's criterion), let $\pi:X\to Y$ be a surjective morphism to a projective manifold. Then $R^i\pi_*(\II(\gamma\varphi))=0$ for $i\geq 1$.
\end{corollary}
\begin{remark} Moreover, by Kodaira's Vanishing Theorem, it is also true that 
$R^i\pi_*(\OO_X)=0$ $(i\geq 1)$ and therefore finally $R^i\pi_*(\OO_{V_\gamma})=0$ $(i\geq 1)$. However, we won't be needing these two properties in the cases we consider.
\end{remark}
We will also need the following simple criterion for the connectedness of fibers of  morphisms. 
\begin{definition}
A morphism between topological spaces is called {\it connected} if all its fibers are connected.
\end{definition}
\begin{proposition}\label{conn_crit}
Let $f:M\to N$ be a proper morphism of complex manifolds. If the morphism $\OO_N\to f_*\OO_M$ is surjective, then $f$ is connected.
\end{proposition}

Based on Proposition \ref{conn_crit}, we now have

\begin{theorem}[{\cite[Theorem 3.1]{Nadel_Wuppertal}}]\label{conn_fibers}
In the situation of Theorem \ref{Nadel_crit} (Nadel's criterion), let $f:X\to Y$ be a connected morphism to a complex projective manifold $Y$. Then 
$(V_\gamma)_y=f^{-1}(y)\cap V_\gamma$ is connected.
\end{theorem}

In Section \ref{678}, we will see that these statement yield Theorem \ref{Nadel_low_degree}. 
\section{The case of 3 points blown up}
Let $X$ be a del Pezzo surface obtained from blowing up three points. W.l.o.g., we assume these to be $P_1=[1,0,0]$, $P_2=[0,1,0]$, $P_3=[0,0,1]$. It is remarked in the introduction to \cite{Nadel_Annals} that this case can be handled easily based on the fact that $X$ is toric. While this is true, we would like to treat this case in line with the (non-toric) cases $r=4,5$, using a finite group of automorphisms.\par
First, we determine the full automorphism group of $X$.\par
Let $H$ be the subgroup of elements of $\PGL(3,\CC)$ that preserve the set $\{P_1,P_2,P_3\}$:
\begin{equation*}
H=\{A\in \PGL(3,\CC):\{A(P_1),A(P_2),A(P_3)\}=\{P_1,P_2,P_3\}\}.
\end{equation*}
The elements of $H$ clearly extend to automorphisms of $X$.\par
Next, consider the standard quadratic Cremona transformation on $\PP^2$:
\begin{equation*}
\Cr([X_0,X_1,X_2])=[X_1X_2,X_0X_2,X_0X_1].
\end{equation*}
$\Cr$ is a birational involution of $\PP^2$ that is a morphism outside the set $\{P_1,P_2,P_3\}$. It lifts to an automorphism of $X$ (see \cite[Example V.4.2.3]{Hartshorne}). We will abuse notation and denote the induced automorphism of $X$ again by $\Cr$.\par
It is a well-known fact (see eg \cite{Koitabashi}) that 
\begin{equation*}
\Aut(X)=H\rtimes \{1,\Cr\}=H\rtimes \ZZ_2.
\end{equation*}
This identity can also be written as
\begin{equation*}
\Aut(X)=(\CC^*\times\CC^*)\rtimes (S_3\times \ZZ_2),
\end{equation*}
where $S_3\subset \PGL(3,\CC)$ is the subgroup of projectivities that act on $\PP^2$ by permuting the coordinates, ie for $\sigma \in S_3$,
\begin{equation*}
\sigma([X_0,X_1,X_2]=[X_{\sigma(0)},X_{\sigma(1)},X_{\sigma(2)}]).
\end{equation*}
Moreover, $\CC^*\times\CC^*$ acts by multiplicitation on, say, the last two coordinates.\par
In order to obtain a finite subgroup of $\Aut(X)$ to use in Criterion \ref{simple_crit}, we note that the Klein four-group $\ZZ_2\times\ZZ_2\subset \CC^*\times \CC^*$ acts on $X$ as follows:
\begin{equation*}
(i,j)([X_0,X_1,X_2])=[X_0,(-1)^iX_1,(-1)^jX_2].
\end{equation*}
Let 
\begin{equation*}
G=(\ZZ_2\times\ZZ_2)\rtimes (S_3\times \ZZ_2).
\end{equation*}
\begin{lemma}\label{r=3_fpf}
The action of $G$ is fixed point free.
\end{lemma}
\begin{proof}
The only fixed point of the action of $S_3$ is $[1,1,1]$. However, this point is not fixed under the action of $(0,1)\in\ZZ_2\times\ZZ_2$. 
\end{proof}
\begin{remark}
It is well-known that $(\ZZ_2\times\ZZ_2)\rtimes S_3=S_4$. From a group theoretic point of view, the statement of Lemma \ref{r=3_fpf} is implied by the fact that $S_4$ admits no complex two-dimensional faithful representations (on the tangent space of a point fof $X$). The non-existence of such representations is a well-known fact in the representation theory of finite groups (see \cite[p. 148]{Dornhoff}).
\end{remark}
\begin{lemma}\label{r=3_eff_invar}
The group $G$ acts effectively on any $G$-invariant irreducible curve.
\end{lemma}
\begin{proof}
For any given element of $G$, it is easy to list the irreducible curves that are left point-wise fixed by the given element (if any exist). However, none of these curves are $G$-invariant.
\end{proof}
We finally observe that $G$ is a group of order $48$ which is clearly not isomorphic to any of the groups listed in Criterion \ref{simple_crit} as finite subgroups of $\Aut(\PP^1)$. Thus, Criterion \ref{simple_crit} can be applied, and we have established the existence of a \ke metric on $X$.
\begin{remark}
Note that the use of the automorphism induced by the quadratic Cremona transformation is crucial. The group $(\ZZ_2\times\ZZ_2)\rtimes S_3=S_4$ is a subgroup of $\Aut(\PP^1)$, and it cannot be used in applying Criterion \ref{simple_crit}.
\end{remark}

\section{The case of 4 points blown up}
Let $X$ be a del Pezzo surface obtained from blowing up four points. We may, and do, assume these to be $P_1=[1,0,0]$, $P_2=[0,1,0]$, $P_3=[0,0,1]$, $P_4=[1,1,1]$. It is known that $\Aut(X)$ is the Weyl group of the root system of Dynkin type $D_5$, which is $S_5$ (see \cite{Koitabashi}). However, we would like to understand $\Aut(X)$ more concretely.  \par
First all of, there is a subgroup $S_4$ of projectivities in $\PGL(3)=\Aut(\PP^2)$ that preserve the set $\{P_1,P_2,P_3,P_4\}$. These projectivities lift to $X$, and we can write $S_4\subset \Aut(X)$.\par
In addition, there exists for every $i=1,\ldots,4$ a quadratic Cremona transformation $\Cr_i$ that leaves $P_i$ fixed and has the three remaining points as indeterminacy locus. Note that such a $\Cr_i$ is only defined up to the action of the $S_3\subset S_4$ consisting of automorphisms fixing $P_i$. For our purposes, it does not matter which $\Cr_i$ we choose. All $\Cr_i$ extend to automorphisms of $X$. In light of this, we can write $\Aut(X)$ set-theoretically as a disjoint union
\begin{equation*}
\Aut(X)=S_4\uplus\left(\biguplus_{i=1}^4  \Cr_i \circ S_4\right).
\end{equation*}
\begin{lemma}\label{r=4_fpf}
The action of $\Aut(X)$ is fixed point free.
\end{lemma}
\begin{proof}
The element of $\Aut(X)$ induced by
\begin{equation*}
\begin{pmatrix}
0&0&1\\
-1&0&1\\
0&-1&1\\
\end{pmatrix}
\end{equation*}
has precisely three fixed points on $X$, namely those points on $X$ lying over 
\begin{equation*}
[1,0,1], [1,1+i,i], [1,1-i,-i] 
\end{equation*}
in $\PP^2$.
On the other hand, the fixed point set of the element of $\Aut(X)$ induced by $\Cr_4$ consists of the three points lying over  
\begin{equation*}
[-1,1,1],[1,1,-1],[1,-1,1]
\end{equation*}
and the $(-1)$-curve over $[1,1,1]$.
Since there is no point appearing in both fixed point sets, the Lemma is proven.
\end{proof}
\begin{remark}
We already noted that $S_4$ admits no complex two-dimensional faithful representations. Since representations clearly remain faithful under restriction to subgroups, $S_5$ does not admit any such representations either. So Lemma \ref{r=4_fpf} follows directly from representation theory, without any explicit computations.
\end{remark} 
\begin{lemma}\label{r=4_eff_invar}
$\Aut(X)$ acts effectively on any $\Aut(X)$-invariant irreducible curve.
\end{lemma}
\begin{proof}
For any given element of $\Aut(X)$, it is easy to list the irreducible curves that are left point-wise fixed by the given element (if any exist). However, none of these curves are $\Aut(X)$-invariant. 
\end{proof}
Lemmas \ref{r=4_fpf} and \ref{r=4_eff_invar} allow us to apply Criterion \ref{simple_crit} to conclude the existence of a \ke metric on $X$.

\section{The case of 5 points blown up}
Let $X$ be a del Pezzo surface obtained by blowing up five points. We can find an automorphism of $\PP^2$ that takes the five points to $P_1=[1,0,0]$, $P_2=[0,1,0]$, $P_3=[0,0,1]$, $P_4=[1,1,1]$, $P_5=[a,b,c]$, with $(a,b,c) \in(\CC^*\times\CC^*\times\CC^*)\backslash\{(1,1,1)\}$. (The reason for $a,b,c \not= 0$ is that no three of these points lie on a line.)\par
The structure of $\Aut(X)$ is worked out in \cite{Hosoh_quartic}. It turns out that it is always of the form
\begin{equation*}
\Aut(X)=\ZZ_2^4\rtimes G_{P_5},
\end{equation*}
where $G_{P_5}$ is a subgroup of $S_5$. The possibilities for $G_{P_5}$ are
\begin{enumerate}
\item $\{\id\},$
\item $\ZZ_2,$
\item $\ZZ_4,$
\item $\ZZ_3\rtimes \ZZ_2,$
\item $\ZZ_5\rtimes \ZZ_2.$
\end{enumerate}
The elements of $G_{P_5}$ are lifts of those elements of $\PGL(3,\CC)$ that preserve the set $\{P_1,P_2,P_3,P_4,P_5\}$. For a generic point $P_5$, we have $G_{P_5}=\{\id\}$. \par
For our purposes, it suffices to argue with elements of $\ZZ_2^4\subseteq \Aut(X)$, regardless of the nature of $G_{P_5}$.\par
\begin{lemma}\label{r=5_fpf}
The action of any subgroup of $\Aut(X)$ isomorphic to $\ZZ_2^3$ or $\ZZ_2^4$  is fixed point free. 
\end{lemma}
\begin{proof}
Instead of using explicit computations, we will argue purely with methods from representation theory. Let $H$ be a subgroup of $\Aut(X)$ isomorphic to either $\ZZ_2^3$ or $\ZZ_2^4$. If $H$ has a fixed point, then it has a {\it faithful} two-dimensional complex representation on the tangent space of the fixed point. \par
Since $H$ is abelian, by Schur's Lemma, the only {\it irreducible} representations of $H$ are maps $\rho : H\to \CC^*$ (see \cite[p.\ 8]{Fulton_Harris}). Since every nontrivial element of $H$ has order two, the image of $\rho$ is either $\{1\}$ or $\{1,-1\}$. Therefore, the image of $H$ in $\GL(2,\CC)$  is isomorphic to either $\{1\}$, $\ZZ_2$, or $\ZZ_2^2$. We have obtained a contradiction.
\end{proof}
Next, let us have a closer look at the elements of $\ZZ_2^4\subseteq \Aut(X)$. The following two birational involutions of $\PP^2$ lift to elements of $\Aut(X)$.\par
We define $\Cr_{45}$  to be
\begin{equation*}
\Cr_{45}([X_0,X_1,X_2])=[aX_1X_2,bX_0X_2,cX_0X_1].
\end{equation*}
This is a quadratic Cremona transformation that exchanges $P_4$ and $P_5$ and has $\{P_1,P_2,P_3\}$ as indeterminacy locus. \par
Moreover, we let $\sigma_1$ be the following cubic involution.
\begin{eqnarray*}
\sigma_1([X_0,X_1,X_2])&=&[-aX_1X_2((c-b)X_0+(a-c)X_1+(b-a)X_2),\\
&&X_1(a(c-b)X_1X_2+b(a-c)X_0X_2+c(b-a)X_0X_1),\\
&&X_2(a(c-b)X_1X_2+b(a-c)X_0X_2+c(b-a)X_0X_1)].
\end{eqnarray*}
The explicit description of $\sigma_1$ tells us that the analogue of Lemmas \ref{r=3_eff_invar} and \ref{r=4_eff_invar} is not true in the case of five points blown up. Namely, it is easy to see that the strict transform of the cubic curve $C$ given by
\begin{equation*}
b(a-c)X_0^2X_2+c(b-a)X_0^2X_1+a(a-c)X_1^2X_2+a(b-a)X_1X_2^2+2a(c-b)X_0X_1X_2=0
\end{equation*}
is precisely the set of points fixed point-wise by the lift of the above $\sigma_1$. In addition, the following Lemma tells us that $C$ is invariant under every element of $\ZZ_2^4$.
\begin{lemma}
Let $A$ be an abelian group acting on a set $M$. For $g\in A$, let 
\begin{equation*}
M^g=\{ x\in M: gx=x\}.
\end{equation*}
Then $M^g$ is invariant under every element of $A$.
\end{lemma}
\begin{proof}
For $x\in M^g$, we have
\begin{equation*}
g(hx)=h(gx)=hx.  
\end{equation*}
\end{proof}
However, there are several reasons why this observation does not impede our proof. \par

Closer inspection (\cite{Koitabashi}, \cite{Hosoh_quartic}) shows that $\ZZ_2^4\subseteq\Aut(X)$ contains the lifts of ten quadratic Cremona involutions $\Cr_{ij}\ (1\leq i<j\leq 5)$ that exchange $P_i,P_j$ and have the remaining three points as indeterminacy locus.  Again, we abuse notation and denote the maps before and after the lift by the same symbols $\Cr_{ij}$.\par

Moreover, $\ZZ_2^4\subseteq\Aut(X)$ contains the lifts of five cubic involutions $\sigma_i\ (1\leq i\leq 5)$. By comparing the respective action on the set of the 16 $(-1)$-curves (which determines any automorphism uniquely, see \cite{Koitabashi} or \cite{Hosoh_quartic}), it is easily verified that 
\begin{eqnarray*}
\Cr_{ij}\circ\Cr_{kl}&=&\Cr_{kl} \circ\Cr_{ij},\\
\Cr_{ij}\circ\Cr_{jk}&=&\Cr_{ik}.
\end{eqnarray*}
Moreover, for $i,j,k,l,m$ all distinct, we have
\begin{equation*}
\sigma_m=\Cr_{ij}\circ\Cr_{kl}.
\end{equation*}
In particular, we have
\begin{equation}\label{group_law}
\sigma_j=\Cr_{1j}\circ\sigma_1 \quad (2\leq j\leq 5).
\end{equation}
\par
The $\Cr_{ij}$ have precisely four fixed points each (both before and after the lifting). Thus, the only $\ZZ_2^4$-invariant curves on which $\ZZ_2^4$ does not act effectively are the lifts of the curves $C$ pertaining to the cubic involutions. However, it is easy to see that for all values $a,b,c\not =0$, the curves $C$ are smooth elliptic curves. Therefore, there exist no smooth rational curves on $X$ on which $G$ acts ineffectively. We omit the easy details of the above argument and instead argue without assuming any knowledge regarding fixed loci of the $\sigma_i$, merely using the group structure of $\ZZ_2^4$.
\begin{lemma}\label{r=5_eff_invar}
The group $\ZZ_2^3$ acts effectively on any $\ZZ_2^4$-invariant irreducible curve.
\end{lemma}
\begin{proof}
Let $C$ be a $\ZZ_2^4$-invariant irreducible curve. Assume that $g\not = \id$ acts trivially on $C$. Then $g=\sigma_i$ for some $1\leq i\leq 5$, because the Cremona maps only have four fixed points each. W.l.o.g., let $i=1$.\par
Because of \eqref{group_law}, none of the $\sigma_j\ (2\leq j\leq 5)$ act trivially on $C$. Thus the group
\begin{equation*}
\ZZ_2^4/\{\id,\sigma_1\} = \ZZ_2^3
\end{equation*}
acts effectively on $C$.
\end{proof}
Lemmas \ref{r=5_fpf} and \ref{r=5_eff_invar} allow us to apply Criterion \ref{simple_crit} to conclude the existence of a \ke metric on $X$.\par
At this point, we have established
\begin{maintheorem12}  [\cite{Siu_ke_symm}, \cite{Tian_1987}, \cite{Tian_Yau}, \cite{Tian_Inv_Math}] Let $X$ be a del Pezzo surface obtained by blowing up $\PP^2$ in $3,4,$ or $5$ points. Then $X$ carries a \ke metric.
\end{maintheorem12}

\section{The cases of 6,7, or 8 points}\label{678}
The following is a concise presentation of \cite[Section 4.1]{Nadel_Wuppertal}.  We have added the missing details for the unspecified ``careful choice" in the proof of \cite[Proposition 4.1]{Nadel_Wuppertal}.

\begin{lemma}
Let $E_1,\ldots,E_r$ be any set of $r$ pairwise disjoint $(-1)$-curves on a del Pezzo surface $X$ of degree $9-r$. Then blowing down $E_1,\ldots,E_r$ gives a birational morphism $X\to\PP^2$.
\end{lemma}
\begin{proof}
The image is a del Pezzo surface of degree $9$. According to the classification, the only such del Pezzo surface is $\PP^2$.
\end{proof}

\begin{lemma}\label{deg_at_least_2}
Let $X$ be a del Pezzo surface of degree $d\leq 4$. Let $C$ be any irreducible curve on $X$. Then there are mutually disjoint $(-1)$-curves $\tilde E_1,\ldots,\tilde E_r$ ($r=9-d$) whose blowing down gives $\tilde \pi:X\to\PP^2$ with $\dim \tilde\pi(C)=1$ and $\deg \tilde\pi(C)\geq 2$.
\end{lemma}
\begin{proof}
Let $E_1,\ldots,E_r$ be the original exceptional curves for the blowing up $\pi:X\to\PP^2$. If $C$ happens to be one of them, say $C=E_1$, then we can contract $E_{12},E_{13},E_{23}$ instead of $E_1=C,E_2,E_3$, where $E_{ij}$ is the strict transform of the line through $P_i$ and $P_j$. Therefore, we can assume w.l.o.g. that $C$ does not get contracted by $\pi$, ie $\dim \pi(C)=1$.\par
If $\pi(C)$ is a line in $\PP^2$ (we are done otherwise), then under any blowing down map, at most two image points of exceptional curves can lie on $\pi(C)$, because the image points will be in general position. In any case, we can assume w.l.o.g. that the image points of $E_3,E_4,E_5$ are each not on $\pi(C)$ (here we use $r\geq 5$). Now let  
\begin{equation*}
\tilde E_1=E_1, \tilde E_2=E_2, \tilde E_3=E_{34}, \tilde E_4=E_{35}, \tilde E_5=E_{45}.
\end{equation*}
If $r\geq 6$, also let
\begin{equation*}
\tilde E_6=E_6,\ldots, \tilde E_r=E_r.
\end{equation*}
We have
\begin{equation*}
\tilde\pi (E_{34}),\tilde\pi(E_{35}), \tilde\pi(E_{45})\in \tilde\pi(C).
\end{equation*}\par
This implies $\deg\tilde\pi(C)\geq 2$, because $\tilde\pi(E_{34}),\tilde\pi(E_{35}), \tilde\pi(E_{45})$ are three points in general position.
\end{proof}
\begin{proposition}\label{no_one_dim_mis}
Let $X$ be a del Pezzo surface of degree $d\leq 3$. Then $-K_X$ does not carry a singular hermitian metric $h$ as in Theorem \ref{Nadel_crit} (Nadel's criterion) such that $\dim V(\II(h))=1$.
\end{proposition}
\begin{proof}
We argue by contradiction. W.l.o.g, let $C=V(\II(h))$ be an irreducible curve. Let $\pi:X\to\PP^2$ be such that $\deg \pi(C)\geq 2$ ($\pi$ exists by Lemma \ref{deg_at_least_2}). For $1\leq i\leq r=9-d$, let $\pi(E_i)=P_i$.\par
{\bf Case 1:} $\deg \pi(C)= 2${\bf .} Since $r\geq 6$, there is $1\leq i\leq r$ such that $P_i\not\in \pi(C)$. W.l.o.g., let $i=r$. Note that $\pi$ factors as
\begin{equation*}
X \xrightarrow{\pi_{1,\ldots,r-1}} \tilde \PP^2 \xrightarrow{\pi_r}\PP^2.
\end{equation*}
The key point is that a general line through $P_r$ intersects $\pi(C)$ (precisely) twice in $\PP^2\backslash\{P_r\}$. Let $\tau$ denote the fibration $\tilde\PP^2\to\PP^1$ whose fibers are lines through $P_r$. Let $f$ denote the composition 
\begin{equation*}
X \xrightarrow{\pi_{1,\ldots,r-1}} \tilde \PP^2 \xrightarrow{\tau}\PP^1.
\end{equation*}
Note that the sets $C_y=f^{-1}(y)\cap C$ have cardinality at least two for general $y\in \PP^1$, which is impossible by Theorem \ref{conn_fibers}.\par

{\bf Case 2:} $\deg \pi(C)\geq 3${\bf .} If there exists an index $1\leq i\leq r$ such that $P_i\not\in \pi(C)$, then we can argue exactly as we did in the previous case. So let us assume that for all $1\leq i\leq r$, $P_i\in \pi(C)$.\par
In order to argue again with the connectedness of fibers, we need to find $i_0$ such that a general line through $P_{i_0}$ intersects $\pi(C)$ in at least two additional distinct points. Let $d=\deg \pi(C)$. If $\mult_{P_{i_0}} \pi(C)\leq d-2$, then the $i_0$ is obviously of the kind we are looking for. We claim that such an $i_0$ always exists. To argue by contradiction, assume that for all $1\leq i\leq r$, we have $\mult_{P_{i_0}} \pi(C)\geq d-1$.\par
Choose a smooth cubic curve $F$ through $P_1,\ldots,P_r$. Then the intersection number of $F$ and $\pi(C)$ is 
\begin{equation*}
F\cdot\pi(C)=3d.
\end{equation*}
It is a basic fact from intersection theory that 
\begin{equation*}
3d=F\cdot\pi(C) \geq \sum_{x\in F\cap\pi(C)} \mult_x{F}\cdot\mult_x \pi(C)\geq\sum_{x\in F\cap\pi(C)} \mult_x\pi(C)\geq r(d-1). 
\end{equation*}
Consequently, $r \geq (r-3)d$. It is easy to check that for $r\geq 6$, this implies $d\leq 2$, a contradiction.
\end{proof}
\begin{theorem13}[\cite{Nadel_Wuppertal}] Let $X$ be a del Pezzo surface of degree no more than three. If $\Aut(X)$ acts without fixed points on $X$, then $X$ has a K\"ahler Einstein metric.
\end{theorem13} 
\begin{proof}
Since there are no zero- or one-dimensional zero-sets of $\Aut(X)$-invariant multiplier ideal sheaves, the theorem follows from Theorem \ref{Nadel_crit} (Nadel's criterion). 
\end{proof}
To clarify the extent to which Theorem \ref{Nadel_low_degree} can actually be applied to prove the existence of \ke metrics on del Pezzo surfaces of degree no more than three, we make the following remarks. An excellent reference with a lucid exposition is \cite[Section 10]{Dolgachev}. \par

Let $X$ be a del Pezzo surface of degree one. The unique base point of the linear system $|-K_X|$ is fixed by all automorphisms. Therefore, unfortunately, $\Aut(X)$ acts with a fixed point regardless of the nature of $X$.\par

On a del Pezzo surface $X$ of degree two, the linear system $|-K_X|$ gives a two-sheeted cover of $\PP^2$ branched along a smooth curve $C$ of degree $4$ in $\PP^2$. This cover defines an involutive automorphism of $X$ called the Geiser involution, which lies in the center of $\Aut(X)$. The corresponding quotient group $G'=\Aut(X)/\ZZ_2$ acts effectively on $C$. If the action of $\Aut(X)$ on $X$ has a fixed point, then so does the action of $G'$ on $C$. At the fixed point, $G'$ has a faithful representation on the one-dimensional tangent space to $C$, so $G'$ must be a cyclic group $\ZZ_n$. Therefore, the action of $\Aut(X)$ has a fixed point only if $\Aut(X)$ is a cyclic central extension of $\ZZ_n$ by $\ZZ_2$. For a list of surfaces $X$ for which this is not the case, see \cite[Table 10.4]{Dolgachev}.\par

For a generic del Pezzo surface $X$ of degree three, $\Aut(X)$ is the trivial group (see \cite{Koitabashi}). However, since del Pezzo surfaces of degree $3$ are precisely the smooth cubic hypersurfaces of $\PP^3$, there are many non-generic $X$ which have extra automorphisms. The following paragraphs give a sufficient criterion for the action of $\Aut(X)$ to be fixed point free.\par
Assume that $x\in X$ is a fixed point of the action of $\Aut(X)$. Then the canonical map $\GL(2,\CC)\to \PGL(2,\CC)$ induces a map $\Aut(X)\to G'=\Aut(X)/Z(\Aut(X))$, where $Z(\Aut(X))$ denotes the center. The tangent plane to $X$ at $x$ intersects $X$ in a singular cubic plane curve $C$, and there are one, two or three lines through $x$ in the tangent plane that are tangent to $C$ at $x$. Let $T\subset \PP^1$ be the set of points corresponding to these lines in the projectivized tangent plane at $x$. Note that $T$ is invariant under the action of $G'$ on $\PP^1$.\par
If $T=\{p_1\}$, then it follows from the classification of finite subgroups of $\PGL(2,\CC)$ mentioned in Section \ref{Nadel_method} that $G'$ is a cyclic group $\ZZ_n$.\par
If $T=\{p_1,p_2\}$, then the classification yields that either $G'$ is a cyclic group $\ZZ_n$ or a dihedral group $\ZZ_n\rtimes \ZZ_2$. \par
If $T=\{p_1,p_2,p_3\}$, then $G'\subset S_3$, ie one of the following holds: $G'=\{\id\}, $ $G'=\ZZ_2,$ $G'=\ZZ_3,$ or $G'=\ZZ_3\rtimes \ZZ_2=S_3$.\par
We have just shown that if $\Aut(X)$ acts with a fixed point, then it is a central extension of one the above $G'$ by a finite cyclic group. Pairing this knowledge with \cite[Table 10.3]{Dolgachev}, it is a simple task to determine which automorphism groups act without fixed points and therefore carry a \ke metric based on Theorem \ref{Nadel_low_degree}. Most notably, the list contains the two examples given by
\begin{equation*}
Z_0^3+Z_1^3+Z_2^3+Z_3^3=0
\end{equation*}
and
\begin{equation*}
Z_0^2 Z_1+Z_2^2 Z_0+Z_3^2 Z_2+Z_1^2 Z_3=0.
\end{equation*}
The automorphism groups are $\ZZ_3^3\rtimes S_4$ resp.~$S_5$. In both cases, the action clearly is fixed point free.\par

\end{document}